\input amstex
\documentstyle{amsppt}
\magnification=\magstep1
\hoffset1 true pc
\voffset2 true pc
\hsize36 true pc
\vsize50 true pc

\tolerance=2000
\define\m1{^{-1}}
\define\ov1{\overline}
\def\gp#1{\langle#1\rangle}
\def\ul2#1{\underline{\underline{#1}}}

\def\rad{ \operatorname{rad \,}}

\TagsOnRight

\topmatter
\title
On existing of filtered multiplicative bases in group algebras
\endtitle
\author
Zsolt Balogh
\endauthor
\leftheadtext\nofrills{Zsolt Balogh} \rightheadtext\nofrills { On
existing of filtered multiplicative bases in group algebras}

\abstract We give an explicit list of all $p$-groups $G$ of order
at most $p^4$ or $2^5$ such that the group algebra $KG$ over the
field $K$ of characteristic $p$ has a filtered multiplicative
$K$-basis.
\endabstract

\address
\hskip-\parindent Institute of Mathematics and Informatics
\newline
College of Ny\'\i regyh\'aza
\newline
S\'ost\'oi \'ut 31/b, H-4410 Ny\'\i regyh\'aza
\newline
Hungary
\newline
baloghzs\@nyf.hu
\endaddress

\thanks
\newline Supported by Hungarian National Fund for Scientific
Research (OTKA) grants No.~T037202 and No.~T043034
\endthanks

\endtopmatter

\document
\subhead 1. Introduction
\endsubhead
In \cite{8} Kupish introduced the following definition. Let $A$ be
a finite-dimensional algebra over a field $K$ and $B$ a $K$-basis
of $A$. Suppose that $B$ is a $K$-basis of $A$ with properties:
\itemitem{1.} if $u,v\in B$ then either $uv=0$ or
$uv\in B$;
\itemitem{2.}  $B\cap \rad(A)$ is a $K$-basis for $\rad(A)$, where $\rad(A)$ denotes  the Jacobson radical of
$A$.

Then $B$ is called a {\it  filtered multiplicative $K$-basis} of
$A$.

R.~Bautista, P.~Gabriel, A.~Roiter and L.~Salmeron showed in
\cite{1} that if there are only finitely many isomorphism classes
of indecomposable $A$-modules over an algebraically closed field
$K$, then $A$ has a filtered multiplicative $K$-basis.

In the present article we shall investigate the following question
from \cite{1}: {\it When have the group algebras $KG$ got a
filtered multiplicative $K$-basis}?

According to Higman's theorem  the group algebra $KG$ over a field
of characteristic $p$ has only finitely many isomorphism classes
of indecomposable $KG$-modules if  and only if all the Sylow
$p$-subgroups of $G$ are cyclic.

Let $G=\gp{a_1}\times \gp{a_2} \times \cdots \times \gp{a_s}$ be a
finite abelian $p$-group with  factors $\gp{a_i}$ of order $q_i$.
Then the set
$$
B=\{(a_1-1)^{n_1}(a_2-1)^{n_2}\cdots (a_s-1)^{n_s}\,\,\,\,
\mid\,\,\,\,   0\leq n_i< q_i\}
$$
forms a filtered multiplicative $K$-basis of the group algebra
$KG$ over the field $K$ of characteristic $p$.

Evidently, if $B_1$ and $B_2$ are filtered multiplicative
$K$-bases of $KG_1$ and $KG_2$, respectively, then $B_1\times B_2$
is a filtered multiplicative $K$-basis of the group algebra
$K[G_1\times G_2]$.

First L.~Paris gave examples of  nonabelian metacyclic $p$-groups
$G$ such that  group algebras  $KG$  have a filtered
multiplicative $K$-bases in \cite{9}.

In \cite{10} P.~Landrock and G.O.~Michler proved that the group
algebra of the smallest Janko group over a field of characteristic
$2$ does not have a  filtered multiplicative $K$-basis.

In \cite{2} the following theorem was proved:

\proclaim {Theorem} Let $G$ be a finite metacyclic $p$-group and
$K$ a field of characteristic $p$. Then the group algebra $KG$
possesses a filtered multiplicative $K$-basis if and only if $p=2$
and exactly one of the following conditions holds:
\itemitem{1.} $G$ is a dihedral group;
\itemitem{2.} $K$ contains a primitive cube root of the unity and
$G$ is a quaternion group of order $8$.
\endproclaim

In \cite{3} was given all $p$-groups $G$ with a cyclic subgroup of
index $p^2$ such that the group algebra $KG$ over the field $K$ of
characteristic $p$ has a filtered multiplicative $K$-basis.

For this question negative answer was given in \cite{3}, when $G$
is either a powerful $p$-group or a two generated $p$-group (
$p\not=2$ ) with  central cyclic commutator subgroup.

\subhead 2. Main results
\endsubhead

Denote $C_{n}$ the cyclic group of order $n$. For the sake of
convenience we shall keep the indices of these groups as in GAP.
We have obtained the following theorems:

\proclaim {Theorem 1} Let $KG$ be the group algebra of a finite
nonabel $p$-group $G$ of order $p^n$  over the field $K$ of
characteristic $p$, where $n< 5$. Then  $KG$ possesses a filtered
multiplicative $K$-basis if and only if $p=2$ and one of the
following conditions satisfy:
\itemitem{1.} $G$ is either dihedral group $D_8$ of order $8$ or dihedral group $D_{16}$ of order $16$;
\itemitem{2.} $G$ is either $Q_8$ or $Q_8\times C_2$  and  $K$  contains
      a primitive cube root of the unity;
\itemitem{3.} $G$ is either $D_8\times C_2$, or the central product $D_8YC_4$ of $D_8$ with $C_4$;
\itemitem{4.} $G$ is $H_{16}=\gp{\quad a,c \quad \vert
\quad a^4=b^2=c^2=1,\;(a,b)=1,\;(a,c)=b,\;(b,c)=1\quad }$.
\endproclaim

\proclaim {Theorem 2} Let $K$ be a field of characteristic $2$ and
 $$
 G=\gp{\quad a,b \quad \vert \quad a^{2^n}=b^{2^m}=c^2=1,\,\, (a,b)=c,\,\, (a,c)=1,\,\,
 (b,c)=1\quad},
 $$
with $n,m\geq 2$. Then  $KG$ possesses a filtered multiplicative
$K$-basis.
\endproclaim

\proclaim {Theorem 3} Let $G$ be the group
$$
\split
G=\gp{\quad a,b \quad  \vert \quad  a^{2^n}=b^2=c^2=d^2=1,\; &(a,b)=c,\;(a,c)=d,\\
                                      & (a,d)=(b,c)=(b,d)=(c,d)=1\quad },\\
\endsplit
$$
with $n>1$, and $K$ a field of characteristic $2$. Then  $KG$ has
no filtered multiplicative $K$-basis.
\endproclaim

\proclaim {Theorem 4} Let $KG$ be the group algebra of a finite
nonabel $2$-group $G$ of order $2^5$  over a field $K$ of
characteristic $2$. Then  $KG$ possesses a filtered multiplicative
$K$-basis if and only if one of the following conditions satisfy:
    \itemitem{1.} $G$ is $G_{18}=D_{32}$, \;  $G_{25}=D_8\times C_4$, \;
        $G_{39}=D_{16}\times C_2$ or $G_{46}=D_{8}\times C_2\times C_2$;
    \itemitem{2.} $G$ is $G_{26}=Q_8\times C_4$, or $G_{47}=Q_8\times C_2\times C_2$  and  $K$  contains
      a primitive cube root of the unity;
    \itemitem{3.} $G$ is  $G_{22}=H_{16} \times C_2$, \; $G_{48}=(D_8YC_4) \times C_2$
    \itemitem{4.} $G$ is one  of the following groups:
\endproclaim
$$
\split
 G_2   =& \gp{\quad a,b \quad \vert \quad a^4=b^4=c^2=1,\,\, (a,b)=c,\,\, (a,c)=1,\,\, (b,c)=1\quad};\\
 G_5   =& \gp{\quad a,b \quad\,  \vert  \quad a^8=b^2=c^2=1,\,\, (a,b)=c,\,\, (a,c)=(b,c)=1\quad};\\
 G_7   =& \gp{\quad a,b,c \,\, \vert \quad a^8=b^2=c^2=1,\,\,(a,c)=a^4,\,\,(a,b)=a^4c,\,\,(b,c)=1\quad };\\
 G_8   =& \gp{\quad a,b,c \,\, \vert \quad a^8=c^2=1,\,\,b^2=a^4,\,\,(a,c)=a^4,\,\,(a,b)=a^4c,\,\,(b,c)=1\quad };\\
 G_9   =& \gp{\quad a,b,c \,\, \vert \quad a^2=b^8=c^2=1,\,\,(b,c)=ab^6,\,\,(a,c)=(a,b)=1\quad };\\
 G_{10}=& \gp{\quad a,b,c \,\, \vert \quad a^8=b^4=c^2=1,\,\,a^4=b^2,\,\,(a,b)=a^6c,\,\,(a,c)=(b,c)=1\quad };\\
 G_{11}=& \gp{\quad a,b,c \,\, \vert \quad a^4=b^4=c^2=1,\,\,(b,c)=ab^2,\,\,(a,c)=(a,b)=1\quad };\\
 G_{49}=& \gp{\quad a,b,c,d \,\, \vert \,\, a^4=1,\,\,b^2=c^2=d^2=a^2,\,\,(a,b)=a^2,\,\,(c,d)=a^2,\\
        &\qquad\qquad\qquad\qquad\qquad\qquad\qquad\qquad (a,c)=(a,d)=(b,c)=(b,d)=1\quad}.
\endsplit
$$

\subhead 3. Preliminary remarks and notation
\endsubhead
Assume that $B$ is a filtered multiplicative $K$-basis for a
finite-dimensional $K$-algebra $A$. In the proof of the  main
results we use the following simple properties  of $B$ (see
\cite{2}):
\itemitem{\,(I)\,\,\,\,} $B\cap rad(A)^n$ is a $K$-basis of
$rad(A)^n$ for all $n\geq 1$.
\itemitem{(II)\,\,} if $u,v\in B\setminus rad(A)^k$ and
$u\equiv v \pmod{rad(A)^k}$ then $u=v$.

Recall  that  the {\it Frattini subalgebra}  $\Phi(A)$ of $A$ is
defined as the intersection of all maximal subalgebras of $A$ if
those exist, and as $A$ otherwise. If $A$ is a nilpotent algebra
over a field $K$, then $\Phi(A)=A^2$ by \cite{5}. It implies that
\itemitem{(III)\,} if $B$ is a filtered multiplicative
$K$-basis  of $A$ and if $B\setminus\{1\}\subseteq rad(A)$,
then all elements of $B\setminus rad(A)^2$ are generators of
$A$ over  $K$.

\bigskip
\noindent A $p$-group $G$ is  called  {\it  powerful}, if one of
the following conditions holds:
\itemitem{$1.$} $G$ is a $2$-group and $G/G^4$ is abelian;
\itemitem{$2.$} $G$ is a $p$-group $(p>2)$ and $G/G^p$ is abelian.

Let $G$ be a finite $p$-group. For $a,b\in G$ we define
$a^b=b\m1ab$ and the commutator $(a,b)=a\m1b\m1ab$.  Denote by
$Q_{2^n}$, $D_{2^n}$ and $SD_{2^n}$ the {\it generalized
quaternion group}, the {\it dihedral } and {\it semidihedral }
$2$-group of order $2^n$, respectively, and
$$
MD_{2^n}= \gp{\quad a,b \quad \vert \quad
a^{2^{n-1}}=b^2=1,\;(a,b)=a^{2^{n-2}} \quad }.$$

We define the {\it  Lazard-Jennings series } $ M_i(G)$ of a finite
$p$-group $G$ by induction ( see \cite{6} ). Put $M_1(G)=G$ and $
M_i(G)=\gp{\,\,\, (M_{i-1}(G),G), M^p_{[{\frac{i}{p}}]}(G))\,\,\,
}, $ where
\itemitem{---} $[\frac{i}{p}]$ is the smallest integer not
less than $\frac{i}{p}$;
\itemitem{---} $\big(M_{i-1}(G), G\big)=
               \gp{\,\,\,  (u,v) \,\,\, \mid \,\,\,
                u\in M_{i-1}(G),\,\,\,   v\in \,\,\,  G\,\,
                }$;
\itemitem{---} $M^p_{i}(G)$ is the  subgroup generated by
$p$-powers of the elements of  $M_{i}(G)$.

Evidently,
$$
M_1(G)\supseteq M_2(G)\supseteq \cdots \supseteq M_t(G)=1.
$$

Let $K$ be a field of characteristic $p$. The ideal
$$
A(KG)=\bigg\{ \quad \sum_{g\in G}\alpha_gg\in KG \quad \mid\quad
\sum_{g\in G}{\alpha_g}=0\,\,\,\, \bigg\}
$$
is called the {\it augmentation} ideal of $KG$. Since $G$ is a
finite $p$-group and $K$ is a field of characteristic $p$, $A(KG)$
is nilpotent, and
$$
A(KG)\supset A^2(KG)\supset\cdots\supset A^{s}(KG)\supset
A^{s+1}(KG)=0.
$$
Moreover, $A(KG)$ is the radical of $KG$.

Then the  subgroup $ {\frak D}_n(G)=\big\{\,\,\, g\in G \,\,\,\,
\mid \,\,\,\,  g-1\in A^n(KG)\,\,\, \big\} $ is called the $n$th
{\it dimensional subgroup}  of $KG$.

It is well known that for finite $p$-group $G$, $M_{i}(G)={\frak
D}_i(G)$ for all $i$.
\bigskip

Let  $\Bbb I=\{i\in \Bbb N\mid {\frak D}_i(G)\not={\frak
D}_{i+1}(G)\}$. For $i\in {\Bbb I}$, let $p^{d_i}$ be the order of
the elementary abelian $p$-group
$$
{\frak D}_i(G)/{\frak D}_{i+1}(G)=\prod_{j=1}^{d_i}
\gp{u_{ij}{\frak D}_{i+1}(G)}.
$$
Hence each $g\in G$ can be written uniquely in the form
$$
g=u_{11}^{\alpha_{11}}u_{12}^{\alpha_{12}}\cdots
u_{1d_1}^{\alpha_{1d_1}} u_{21}^{\alpha_{21}} \cdots
u_{2d_2}^{\alpha_{2d_2}}\ldots u_{i1}^{\alpha_{i1}}
\cdots u_{id_i}^{\alpha_{id_i}}\cdots
u_{s1}^{\alpha_{s1}}\cdots
u_{sd_s}^{\alpha_{sd_s}},
$$
where the indices are in lexicographic order, $i\in \Bbb I$,
$0\leq \alpha_{ij}<p$, and $s$ is defined  as above.

Let $w=\prod_{l\in \Bbb
I}(\prod_{k=1}^{d_l}(u_{lk}-1)^{y_{lk}})\in A(KG)$ be where $0\leq
y_{lk}<p$, and the indices of the factors are in lexicographic
order. Then $w$ is called a regular element of weight
$\mu(w)=\sum_{l\in {\Bbb I}}\big(\sum_{k=1}^{d_l}ly_{lk}\big)$. By
Jennings Theorem ( see \cite{6} ), regular elements which weight
not less than $t$ constitute a $K$-basis for the ideal $A^t(KG)$.

\noindent Clearly, $\{\,\,\,(u_{1j}-1)+A^2(KG)\,\,\, \mid\,\,\,
j=1,\ldots,d_1\,\,\, \}$ is a $K$-basis of $A(KG)/A^2(KG)$.

Note that ${\frak D}_2(G)$ coincides with the Frattini subgroup of
$G$, so the set $\{u_{11},u_{12},\ldots,u_{1d_1}\}$ is a minimal
generator system of $G$.

Suppose that $B_1=\{1 \}\cup \{b_1,b_2,\ldots,b_{|G|-1} \}$ is a
filtered multiplicative $K$-basis for $KG$. Then $B=B_1\setminus
\{1\}$ is a filtered multiplicative $K$-basis  of $A(KG)$ and
contains $|G|-1$ elements.

Let $B\setminus(B\cap A^2(KG))=\{b_1,b_2,\ldots,b_n\}$. Evidently,
$n=d_1$ and
$$
b_k\equiv \sum_{i=1}^n\alpha_{ki}(u_{1i}-1)\pmod {A^2(KG)},
$$
where
$\alpha_{ki}\in K$ and $\Delta=det(\alpha_{ki})\not=0$.

For units $x,y$ of $KG$ we have
$$
\split
(y-1)(x-1)&=[(x-1)(y-1)+(x-1)+(y-1)](z-1)\\
           &+(x-1)(y-1)+(z-1),
\endsplit\tag1
$$
where $z=(y,x)$. Since $z_{ji}=(u_{1j},u_{1i})\in {\frak D}_2(G)$
and $z_{ji}-1\in A^2(KG)$, using $(1)$ we obtain that
$$
(u_{1j}-1)(u_{1i}-1)\equiv (u_{1i}-1)(u_{1j}-1)+(z_{ji}-1)
\pmod{A^3(KG)}.\tag2
$$

Thus simple computations give that
$$
\split
b_kb_s\equiv\sum_{i=1}^n\alpha_{ki}\alpha_{si}(u_{1i}-1)^2
&+
\sum_{i,j=1\atop {i<j}}^n(\alpha_{ki}\alpha_{sj}+\alpha_{kj}
\alpha_{si})(u_{1i}-1)(u_{1j}-1)\\
& +\sum_{i,j=1\atop i<j}^n\alpha_{kj}\alpha_{si}(z_{ji}-1)
\pmod{A^3(KG)},
\endsplit
\tag3
$$
where $k,s=1,\ldots,n$.

Denote by $\frak A$ the set of groups which belong to one of the
following type of nonabelian  $p$-groups:
\item{1.} either metacyclic or powerful;
\item{2.} $p$-group with cyclic subgroup of index $p^2$;
\item{3.} two generated $p$-group $(p\ne 2)$ with central cyclic
commutator subgroup.

\subhead
{4. Proof of Theorem 1}
\endsubhead
Let $K$ be a field of characteristic $p$ ($p$ is odd) and $G$ a
$p$-group of order $p^4$. The classification of these groups can
be found in \cite{7}. According to \cite{2} and \cite{3} if $G$
belongs to $\frak{A}$ then $G$ has no filtered multiplicative
basis. If $G$ does not belong to $\frak{A}$, then it is one of the
following two groups:
$$
\split H_1=\gp{\quad a,c \quad \vert  \quad
a^p=c^p=1, \quad  & (a,c)=d,(d,c)=f,\\
&(a,d)=(a,f)=(c,f)=(d,f)=1\quad}\\
\endsplit
\text{ with } p>3;
$$
and
$$
\split
H_2=\gp{\quad a,c \quad  \vert \quad   a^p=c^p=1, \,  & (a,c)=d,\\
                                                        &(c,d)=(a,d)=1\quad} \times
    \gp{\quad h \,\, \vert \,\, h^p=1\quad }\\
\endsplit
\text{ with } p\geq 3.
$$
It is easy to check that in both group algebras $KH_1$ and $KH_2$:
$$
(c-1)(a-1)\equiv (a-1)(c-1)-(d-1) \pmod{A^3(KG)}.\tag4
$$

Let us consider the following cases:

{\bf Case 1.} Let $G=H_1$. Since
$$
M_1(G)=G,\quad  M_2(G)=\gp{G',G^p}=\gp{d,f},\quad M_3(G)=\gp{\;
(\gp{d,f},G),\; G^p\;}=\gp{f}
$$
we have that $\mu(d)=2$ and $\mu(f)=3$, where $\mu$ is the weight
of these elements. Using (4) and
$$
(d-1)(c-1)\equiv (c-1)(d-1)+(f-1) \pmod{A^4(KG)},
$$
  let us compute $b_{i_1}b_{i_2}b_{i_3}$  modulo
${A^4(KG)}$ where ($i_k=1,2$). The results  of our  computations
will be written  in a  table, consisting of the coefficients of
the decomposition $b_{i_1}b_{i_2}b_{i_3}$ with respect to the
basis
$$
\split \{\quad  (a-1)^{j_1}(c-1)^{j_2}(d-1)^{j_3}(f-1)^{j_4} \quad
\mid \quad  &j_1+j_2+2j_3+3j_4=3;\\
& j_1,j_2=0,1,2,3;\; j_3,j_4=0,1\quad \}
\endsplit
$$
of the ideal $A^3(KG)/A^4(KG)$. The coefficients of
$b_{i_1},b_{i_2},b_{i_3}$ will be denoted $\alpha_i, \beta_i,
\gamma_i$, respectively, and in the following we shall use these
coefficients. We shall divide the table into two parts (the second
part written below the first part). The coefficients corresponding
to the first four basis elements will be in the first part of the
table, while the next three will be in the second one. Thus
$$
\eightpoint \vbox {{ \offinterlineskip
 \halign {
$#$ \ &  \vrule   \    $#$ \ &  \vrule   \    $#$ \ & \vrule \ $#$
\ &  \vrule   \    $#$ \  \cr
 &
(a-1)^3 &(a-1)^2(c-1) & (a-1)(d-1) & (a-1)(c-1)^2 \cr \noalign
{\hrule}
 b_1b_2b_1 & \alpha_1^2\beta_1 &
2\alpha_1\alpha_2\beta_1+\alpha_1^2\beta_2 &
-2\alpha_1\alpha_2\beta_1-\alpha_1^2\beta_2 &
2\alpha_1\alpha_2\beta_2+\alpha_2^2\beta_1  \cr
b_1b_2^2 & \alpha_1\beta_1^2 &
2\alpha_1\beta_1\beta_2+\alpha_2\beta_1^2 &
-2\alpha_2\beta_1^2-\alpha_1\beta_1\beta_2 &
2\alpha_2\beta_1\beta_2+\alpha_1\beta_2^2  \cr
b_2b_1^2 & \alpha_1^2\beta_1 &
2\alpha_1\alpha_2\beta_2+\alpha_2^2\beta_1 &
-2\alpha_1^2\beta_2-\alpha_1\alpha_2\beta_1 &
2\alpha_1\alpha_2\beta_2+\alpha_2^2\beta_1  \cr
b_2b_1b_2 & \alpha_1\beta_1^2 &
2\alpha_1\beta_1\beta_2+\alpha_2\beta_1^2 &
-2\alpha_1\beta_1\beta_2-\alpha_2\beta_1^2 &
2\alpha_2\beta_1\beta_2+\alpha_1\beta_2^2  \cr
b_1^3 & \alpha_1^3 & 3\alpha_1^2\alpha_2 & -3\alpha_1^2\alpha_2 &
3\alpha_1\alpha_2^2  \cr
b_1^2b_2 & \alpha_1^2\beta_1 &
2\alpha_1\alpha_2\beta_1+\alpha_1^2\beta_2 &
-3\alpha_1\alpha_2\beta_1 &
2\alpha_1\alpha_2\beta_2+\alpha_2^2\beta_1  \cr
b_2^2b_1 & \alpha_1\beta_1^2 &
2\alpha_1\beta_1\beta_2+\alpha_2\beta_1^2 &
-3\alpha_1\beta_1\beta_2 &
2\alpha_2\beta_1\beta_2+\alpha_1\beta_2^2  \cr
b_2^3 & \beta_1^3 & 3\beta_1^2\beta_2 & -3\beta_1^2\beta_2 &
3\beta_1\beta_2^2  \cr \noalign {\hrule} }}}
$$
$$
\eightpoint \vbox {{ \offinterlineskip
 \halign {
$#$ \ &  \vrule   \    $#$ \ &  \vrule   \    $#$ \ & \vrule \ $#$
\  \cr
 &
(c-1)(d-1)& (f-1) & (c-1)^3 \cr \noalign {\hrule}
 b_1b_2b_1 &
-2\alpha_1\alpha_2\beta_2-\alpha_2^2\beta_1 &
-\alpha_2^2\beta_1-\alpha_1\alpha_2\beta_2 & \alpha_2^2\beta_2 \cr
b_1b_2^2 & -3\alpha_2\beta_1\beta_2 & -2\alpha_2\beta_1\beta_2 &
\alpha_2\beta_2^2 \cr
b_2b_1^2 & -3\alpha_1\alpha_2\beta_2 & -2\alpha_1\alpha_2\beta_2 &
\alpha_2^2\beta_2 \cr
b_2b_1b_2 & -2\alpha_2\beta_1\beta_2 -\alpha_1\beta_2^2 &
-\alpha_1\beta_2^2-\alpha_2\beta_1\beta_2 & \alpha_2\beta_2^2 \cr
b_1^3 &  -3\alpha_1\alpha_2^2 & -2\alpha_1\alpha_2^2 & \alpha_2^3
\cr
b_1^2b_2 & -2\alpha_2^2\beta_1-\alpha_1\alpha_2\beta_2 &
-\alpha_1\alpha_2\beta_2-\alpha_2^2\beta_1 & \alpha_2^2\beta_2 \cr
b_2^2b_1 & -2\alpha_1\beta_2^2-\alpha_2\beta_1\beta_2 &
-\alpha_2\beta_1\beta_2-\alpha_1\beta_2^2 & \alpha_2\beta_2^2 \cr
b_2^3 &  -3\beta_2\beta_2^2 & -2\beta_1\beta_2^2 & \beta_2^3 \cr
\noalign {\hrule} }}}
$$

We have obtained $8$ elements, but the $K$-dimension of
$A^3(KG)/A^4(KG)$ equals $7$. Since $\Delta \neq 0$, we can
establish that one of this elements either equals to zero modulo
ideal $A^4(KG)$ or coincides with another one.

It is easy to see that none of lines are equal to zero. Indeed,
for example, if $b_1b_2b_1 \equiv 0 \pmod{A^4(KG)}$ then from
second column of the first part and fourth column of the second
part of the table we get that $\alpha_1^2\beta_1=0$ and
$\alpha_2^2\beta_2=0$. Since $\Delta \ne 0$, this case is
impossible by third column of the first part of this table. In a
similar manner we can proof this statement for all lines.

The assumption that two of lines are equal also a contradiction.
For instance, if $b_1b_2b_1 \equiv b_1b_2^2 \pmod{A^4(KG)}$, then
from second column of the first part and fourth column of the
second part of the table it follows that
$\alpha_1\beta_1(\alpha_1-\beta_1)=0$ and
$\alpha_2\beta_2(\alpha_2-\beta_2)=0$. Since $\Delta \ne 0$, the
third column of the first part of the table leads to a
contradiction.

Similar calculations for any two lines also lead to a
contradiction, so we have got that $KG$ has no filtered
multiplicative basis.

{\bf Case 2.} Let $G=H_2$. Using (4)  let us compute
$b_{i_1}b_{i_2}$ modulo ${A^3(KG)}$ where ($i_k=1,2,3$). The
results of our  computations will be written  in a  table,
consisting of the coefficients of the decomposition
$b_{i_1}b_{i_2}$ with respect to the basis
$$
\split \{\quad  (a-1)^{j_1}(c-1)^{j_2}(h-1)^{j_3}(d-1)^{j_4} \quad
\mid \quad  &j_1+j_2+j_3+2j_4=2;\\
& j_1,j_2,j_3=0,1,2;\,\, j_4=0,1\quad \}
\endsplit
$$
of the ideal $A^2(KG)/A^3(KG)$:
$$
\eightpoint \vbox {{ \offinterlineskip
 \halign {
\; $#$ & \vrule \; $#$ \; & \vrule \; $#$ \; & \vrule \; $#$ \; &
\vrule \; $#$ \; & \vrule \; $#$ \; & \vrule \; $#$ \; & \vrule \;
$#$ \; & \vrule \; $#$ \; \cr & \scriptstyle{(a-1)^2} &
\scriptstyle{(a-1)(c-1)} & \scriptstyle{(a-1)(h-1)} &
\scriptstyle{(c-1)^2} & \scriptstyle{(c-1)(h-1)} &
\scriptstyle{(h-1)^2} & \scriptstyle{(d-1)}\cr \noalign {\hrule}
b_1b_2                                            &
\alpha_1\beta_1 & \alpha_1\beta_2+\alpha_2\beta_1 &
\alpha_1\beta_3+\alpha_3\beta_1& \alpha_2\beta_2 &
\alpha_2\beta_3+\alpha_3\beta_2& \alpha_3\beta_3 &
-\alpha_2\beta_1 \cr
b_2b_1                                            &
 \alpha_1\beta_1                              &
\alpha_1\beta_2+\alpha_2\beta_1                    &
\alpha_1\beta_3+\alpha_3\beta_1& \alpha_2\beta_2 &
\alpha_2\beta_3+\alpha_3\beta_2 & \alpha_3\beta_3 &
-\alpha_1\beta_2 \cr
b_1b_3                                             &
\gamma_1\alpha_1 & \alpha_1\gamma_2+\alpha_2\gamma_1 &
\alpha_1\gamma_3+\alpha_3\gamma_1                    &
\alpha_2\gamma_2 &   \alpha_2\gamma_3+\alpha_3\gamma_2 &
\alpha_3\gamma_3 & -\alpha_2\gamma_1 \cr
b_3b_1                                            &
\gamma_1\alpha_1 &   \alpha_1\gamma_2+\alpha_2\gamma_1 &
\alpha_1\gamma_3+\alpha_3\gamma_1& \alpha_2\gamma_2 &
\alpha_2\gamma_3+\alpha_3\gamma_2& \alpha_3\gamma_3 &
-\alpha_1\gamma_2 \cr
b_2b_3                                            &
\beta_1\gamma_1 & \beta_1\gamma_2+\beta_2\gamma_1 &
 \beta_1\gamma_3+\beta_3\gamma_1              & \beta_2\gamma_2 &
 \beta_2\gamma_3+\beta_3\gamma_2              & \beta_3\gamma_3 &
-\beta_2\gamma_1 \cr
b_3b_2                                            &
\beta_1\gamma_1 & \beta_1\gamma_2+\beta_2\gamma_1 &
\beta_1\gamma_3+\beta_3\gamma_1 & \beta_2\gamma_2 &
\beta_2\gamma_3+\beta_3\gamma_2 & \beta_3\gamma_3 &
-\beta_1\gamma_2 \cr
b_1^2                                           & \alpha_1^2 &
2\alpha_1\alpha_2                             & 2\alpha_1\alpha_3
& \alpha_2^2 & 2\alpha_2\alpha_3                             &
\alpha_3^2 & -\alpha_1\alpha_2 \cr
b_2^2                                           & \beta_1^2 &
2\beta_1\beta_2                               & 2\beta_1\beta_3 &
\beta_2^2                                     & 2\beta_2\beta_3 &
\beta_3^2                                     & -\beta_1\beta_2
\cr
b_3^2                                           & \gamma_1^2 &
2\gamma_1\gamma_2                             & 2\gamma_1\gamma_3
& \gamma_2^2 & 2\gamma_2\gamma_3                             &
\gamma_3^2 & -\gamma_1\gamma_2 \cr \noalign {\hrule} }}}
$$
We have obtained $9$ elements, but the $K$-dimension of
$A^2(KG)/A^3(KG)$ equals $7$, so we conclude that some lines of
the table either are  equal to zero modulo the ideal $A^3(KG)$ or
coincide with some other lines.

Since $\Delta \neq 0$, it is clear that $b_i^2\not\equiv 0$ and
$b_ib_j\not\equiv 0 \pmod{A^3(KG)}$. According to the last tree
lines of the table if $b_i^2\equiv b_j^2 \pmod{A^3(KG)}$, then
either $b_i\equiv b_j$ or $b_i\equiv -b_j \pmod{A^3(KG)}$, so we
have that $b_1b_2 \equiv b_2b_1$, $b_1b_3 \equiv b_3b_1$ and
$b_2b_3 \not\equiv b_3b_2 \pmod{A^3(KG)}$, because the other cases
are similar to this one.

Simple computations show that if either $\alpha_1\ne 0$ or
$\beta_1\ne 0$, then $KG$ is a commutative algebra which is a
contradiction, so we can assume that $\alpha_1=\beta_1=0$. From
the 8th column we have $\alpha_2\gamma_1=0$. Since $\Delta \neq 0$
we conclude that $\alpha_2=0$ and we have a basis of
$A(KG)/A^2(KG)$:
$$
\cases
 b_1 &\equiv (h-1) \pmod{A^2(KG)};               \\
 b_2 &\equiv (c-1)+\beta_3(h-1) \pmod{A^2(KG)};    \\
 b_3 &\equiv (a-1)+\gamma_2(c-1)+\gamma_3(h-1) \pmod{A^2(KG)}.      \\
 \endcases
$$

Let us compute $b_{i_1}b_{i_2}b_{i_3}$ modulo ${A^4(KG)}$ where
$i_k=1,2,3$ with respect to the basis
$$
\split \{\quad  (a-1)^{j_1}(c-1)^{j_2}(h-1)^{j_3}(d-1)^{j_4} \quad
\mid \quad  &j_1+j_2+j_3+2j_4=3;\\
& j_1,j_2,j_3=0,1,2,3;\,\, j_4=0,1\quad \}
\endsplit
$$
of the ideal $A^3(KG)/A^4(KG)$.

Assume that $p=3$. Since the dimension of $A^3(KG)/A^4(KG)$ is
$10$, so we conclude that
$$
b_1^2b_2 \equiv b_1b_2^2, \quad b_2^2b_3 \equiv b_3b_2^2, \quad
b_3^2b_2 \equiv b_3b_2b_3 \pmod{A^4(KG)},
$$
and $b_1^3 \equiv 0, \quad b_2^3 \equiv 0 \pmod{A^4(KG)}$. From
these congruences give that $\beta_3=\gamma_2=\gamma_3=0$.

Now suppose that $p>3$. In this case the dimension of
$A^3(KG)/A^4(KG)$ is $15$, so we conclude that
$$
b_1^2b_2 \equiv b_1b_2^2, \quad b_3^2b_2 \equiv b_3b_2b_3
\pmod{A^4(KG)},
$$
and we also get that $\beta_3=\gamma_2=\gamma_3=0$.

Assume that $KG$ has a filtered multiplicative basis. Since $KG=K[
G_1 \times G_2 ]$, where
$$
G_1=\gp{\quad a,c \quad
\vert \quad   a^p=c^p=1, \quad (a,c)=d,(c,d)=(a,d)=1\quad}
$$
and $G_2=\gp{\quad h \quad \vert \quad h^p=1\quad}$, and we have
established that
$$
\cases
 b_1 &\equiv (h-1) \pmod{A^2(KG)};               \\
 b_2 &\equiv (c-1) \pmod{A^2(KG)};    \\
 b_3 &\equiv (a-1) \pmod{A^2(KG)},      \\
 \endcases
$$
with  $b_1 \in KG_1$, $b_2,b_3 \in KG_2$, so we conclude that
$KG_2$ also has filtered multiplicative basis, which is a
contradiction by \cite{3}.

Let $K$ be a field of characteristic $2$. If $|G|<2^5$, then $KG$
has a filtered multiplicative basis (see  \cite{2,3}) if and only
if $G$ and $K$ satisfy the conditions of Theorem $1$, so the proof
of the theorem is complete.


\subhead {5. Proof of Theorem 2}
\endsubhead
 Let
 $$
 G=\gp{\quad a,b \quad \vert \quad a^{2^n}=b^{2^m}=c^2=1,\,\, (a,b)=c,\,\, (a,c)=1,\,\,
 (b,c)=1\quad},
 $$
 and put
$$
b_1^1=u \equiv (1+a),\quad b_1^2=v \equiv (1+b) \pmod{A^2(KG)}.
$$
Using the identity:
$$\split
(1+b)(1+a) &\equiv  (1+a)(1+b)+(1+c) \pmod{A^3(KG)},
\endsplit
$$
we get that the set $\{\; b_2^1 =uv, \; b_2^2 =vu,\; b_2^3 =u^2,\;
b_2^4 =v^2\; \}$ is a basis of $A^2(KG)/A^3(KG)$ and
$$
\split b^1_3&=uvu\equiv (1+a)^2(1+b)+(1+a)(1+c) \pmod{A^3(KG)};\\
b^2_3&=u^2v\equiv (1+a)^2(1+b) \pmod{A^3(KG)};\\
b^3_3&=u^3\equiv (1+a)^3 \pmod{A^3(KG)};\\
b^4_3&=uv^2\equiv (1+a)(1+b)^2 \pmod{A^3(KG)};\\
b^5_3&=vuv\equiv (1+a)(1+b)^2+(1+b)(1+c) \pmod{A^3(KG)};\\
b^6_3&=v^3\equiv (1+b)^3 \pmod{A^3(KG)};\\
\endsplit
$$
is a basis for $A^3(KG)/A^4(KG)$ and its determinant $\Delta_3=1$.
We shall construct a basis of $A^i(KG)/A^{i+1}(KG)$ by induction.
Assume that $b^1_{i-1},b^2_{i-1},\cdots,b^{n-1}_{i-1},b^n_{i-1}$
is a basis for $A^{i-1}(KG)/A^i(KG)$. Evidently, the determinant
$\Delta_{i-1}$ of this basis is not zero. Simple computations show
that the determinant $\Delta_i$ of the elements
$b^j_i=ub^j_{i-1}$, for $j=\{1,2,\cdots,n \}$ and
$b^{n+1}_i=b^{n-1}_{i-1}v$, $b^{n+2}_i=b^{n}_{i-1}v$ is
equal to $\Delta_{i-1}\cdot \left| \smallmatrix 1 & 0 \\ 0 & 1 \\
\endsmallmatrix \right| \ne 0$, so we got $n+2$ linearly independent elements. Since $dim\,
A^{i}(KG)/A^{i+1}(KG)$ is also $n+2$ we have obtained that $KG$
has a filtered multiplicative basis.

\subhead {6. Proof of Theorem 3}
\endsubhead
Let $G$ be the group
$$
\split
G_6=\gp{\quad a,b \quad  \vert \quad  a^{2^n}=b^2=c^2=d^2=1,\; &(a,b)=c,\;(a,c)=d,\\
                                      & (a,d)=(b,c)=(b,d)=(c,d)=1\quad },\\
\endsplit
$$
with $n>1$. Let us compute the Lazard-Jennings series of this
group:
$$
\quad M_1(G)=G, \quad M_2(G)=\gp{a^2, c,d},\quad
M_3(G)=\gp{d},\quad M_4(G)=\gp{1}.
$$

We conclude that $\mu(c)=2$ and $\mu(d)=3$. Using the identity
$$
(1+b)(1+a) \equiv (1+a)(1+b)+(1+c) \pmod{A^3(KG)},\tag5
$$
it follows that
$$
\eightpoint
       \cases
&b_1b_2 \equiv \alpha_1\beta_1(1+a)^2+(\alpha_1\beta_2+\alpha_2\beta_1)(1+a)(1+b)+\alpha_2\beta_1(1+c) \pmod{A^3(KG)}; \\
&b_2b_1 \equiv \alpha_1\beta_1(1+a)^2+(\alpha_1\beta_2+\alpha_2\beta_1)(1+a)(1+b)+\alpha_1\beta_2(1+c) \pmod{A^3(KG)}; \\
&b_1^2 \equiv \alpha_1^2(1+a)^2+\alpha_1\alpha_2(1+c) \pmod{A^3(KG)}; \\
&b_2^2 \equiv \beta_1^2(1+a)^2+\beta_1\beta_2(1+c) \pmod{A^3(KG)}.
\endcases
$$
We have obtained $4$ elements, but the $K$-dimension of
$A^2(KG)/A^3(KG)$ equals $3$. Since $\Delta \neq 0$, we get that
$b_1b_2 \not\equiv b_2b_1,b_1^2,b_2^2\;$ and $b_1b_2,\; b_2b_1
\not\equiv 0$ and $b_1^2 \not\equiv b_2^2 \pmod{A^3(KG)}$. Thus
either $b_1^2 \equiv 0$ or $b_2^2 \equiv 0 \pmod{A^3(KG)}$. It is
easy to see that the second case is similar to the first one, so
we consider the second one. Let $\beta_1=0$ and we can put
$\alpha_1=\beta_2=1$ and
$$
\split u=b_1 &\equiv (1+a)+\alpha_2(1+b)\pmod{A^2(KG)};
\\  v=b_2 &\equiv (1+b)\pmod{A^2(KG)}. \\
\endsplit
$$
Using (5) and the identity
$$
(1+c)(1+a) \equiv (1+a)(1+c)+(1+d) \pmod{A^4(KG)},
$$
straightforward computations show that
$$
\eightpoint \cases
uvu^2 &\equiv (1+a)^3(1+b)+\alpha_2(1+a)(1+b)(1+c)+(1+a)(1+d)         \pmod{A^5(KG)};\\
vu^3  &\equiv (1+a)^3(1+b)+(1+a)^2(1+c)+\alpha_2(1+a)(1+b)(1+c)+ \\
      &\qquad \qquad \qquad \qquad  (1+a)(1+d)                                 \pmod{A^5(KG)};\\
vuvu  &\equiv (1+b)(1+d)+(1+a)(1+b)(1+c)                              \pmod{A^5(KG)};\\
u^2vu &\equiv (1+a)^3(1+b)+(1+a)^2(1+c)+\alpha_2(1+a)(1+b)(1+c)+
\\
      &\qquad \qquad \qquad \qquad \alpha_2(1+b)(1+d) \pmod{A^5(KG)};\\
uvuv  &\equiv (1+a)(1+b)(1+c)                                         \pmod{A^5(KG)};\\
vu^2v &\equiv (1+b)(1+d)                                              \pmod{A^5(KG)};\\
u^3v  &\equiv
(1+a)^3(1+b)+\alpha_2(1+a)(1+b)(1+c)+\alpha_2(1+b)(1+d)
\pmod{A^5(KG)}.
\endcases
$$

We have obtained 7 different element, but this is a contradiction
because $dim(A^4(KG)/A^5(KG))=5$.


\subhead {6. Proof of Theorem 4}
\endsubhead

Let $G$ be  a nonabelian $2$-group of order $2^5$. According to
\cite{3} if $G$ is one of the following groups
$\{G_5,G_7,G_8,G_9,G_{10},G_{11} \}$, then $G$ has a cyclic
subgroup of index $p^2$ and $KG$ has filtered multiplicative
basis, but if $G$ is one of the following groups:
$$ \split
 G_{40}   =&SD_{16} \times C_2;\\
 G_{41}   =&Q_{16} \times C_2;\\
 G_{42}   =&\gp{\quad a,b,c \quad \vert a^8=b^4=c^4=1,a^4=b^2=c^2,(a,b)=a^6,(a,c)=(b,c)=1} ;\\
 G_{43}   =&\gp{\quad a,b,c \quad \vert a^8=b^2=c^2=1,(a,b)=a^6,(a,c)=a^4,(b,c)=1} ;\\
 G_{44}   =&\gp{\quad a,b,c \quad \vert a^8=c^2=1,b^2=a^4,(a,c)=a^4,(a,b)=a^6,(b,c)=1},\\
\endsplit
$$
then $KG$ has no filtered multiplicative basis.

If $G$ is one of the following groups:
$$ \split
 G_4      =&\gp{\quad a,b \quad \vert \quad a^8=b^4=1,\,\,(a,b)=a^4\quad};\\
 G_{37}   =&MD_{16} \times C_2;\\
 G_{38}   =&\gp{\quad a,b,c \quad \vert \quad a^8=b^2=c^2=1,\,\, (b,c)=a^4,\,\,(a,b)=(a,c)=1\quad}.\\
\endsplit
$$
then they are powerful groups and by \cite{3} $KG$ has no a
filtered multiplicative basis. If $G$ is one of the following
groups:
$$ \split
 G_{12}   =&\gp{\quad a,b \quad \vert a^4=b^8=1,(a,b)=a^2} ;\\
 G_{13}   =&\gp{\quad a,b \quad \vert a^8=b^4=1,(a,b)=a^2} ;\\
 G_{14}   =&\gp{\quad a,b \quad \vert a^8=b^4=1,(a,b)=a^6} ;\\
 G_{15}   =&\gp{\quad a,b \quad \vert a^8=1,b^4=a^4,(a,b)=a^6} ;\\
 G_{17}=&MD_{32},\;G_{18}=D_{32},\;G_{19}=SD_{32},\;G_{20}=Q_{32},\\
\endsplit
$$
then $G$ is a metacyclic group and $KG$ has a filtered
multiplicative basis if and only if $G=G_{18}$ by \cite{2}.

According to \cite{2} and \cite{3} we get that for the following
direct products $KG$ has a filtered multiplicative basis:
$G_{22}=H_{16} \times C_2, G_{25}=D_8 \times C_4, G_{26}=Q_8\times
C_4, G_{39}=D_{16}\times C_2, G_{46}= D_8\times C_2 \times C_2,
G_{47}= Q_8\times C_2 \times C_2, G_{48}=(D_8YC_4)\times C_2.$


If $G=G_{2}$, then for $n=m=2$  Theorem $2$ asserts that $KG$ has
a filtered multiplicative $K$-basis.

Let $G$ be the group
$$
\split
G_6=\gp{\quad a,b \quad  \vert \quad  a^4=b^2=1,&\; (a,b)=c,\;(a,c)=d,\\
                                      & (a,d)=(b,c)=(b,d)=(c,d)=1\quad }.\\
\endsplit
$$
For $n=2$ the group in Theorem $3$ is isomorphic to $G_6$, so the
group algebra $KG_6$ has no filtered multiplicative $K$-basis.

Now, we shall consider the following $7$ cases.


{\bf Case 1.} Let $G$ be the group
$$
G_{23}=\gp{\quad a,b,c \quad  \vert \quad  a^4=b^4=c^2=1,\;
(a,c)=(b,c)=1,\;(a,b)=a^2\quad}.
$$

Using the identity
$$
(1+b)(1+a) \equiv (1+a)(1+b)+(1+a)^2 \pmod{A^3(KG)},
$$
let us compute $b_{i_1}b_{i_2}$  modulo ${A^3(KG)}$ where
($i_k=1,2,3$). The results of our  computation will be written  in
a  table, consisting of the coefficients of the decomposition
$b_{i_1}b_{i_2}$ with respect to the  basis
$$
\split \{\quad  (1+a)^{j_1}(1+b)^{j_2}(1+c)^{j_3} \quad
\mid \quad  &j_1+j_2+j_3=2;\\
& j_1,j_2=0,1,2;\; j_3=0,1 \quad \}
\endsplit
$$
of the ideal $A^2(KG)$.

$$
\eightpoint \vbox {{ \offinterlineskip
 \halign {
$#$ \; &  \vrule   \; $#$ \; &  \vrule   \; $#$ \; &  \vrule   \;
$#$ \; &  \vrule   \; $#$ \; &  \vrule   \; $#$ \; \cr
    &  (1+a)^2    &(1+a)(1+b) &  (1+a)(1+c) & (1+b)(1+c) & (1+b)^2 \cr
\noalign {\hrule}
b_1b_2                                                    &
\alpha_1\beta_1+\alpha_2\beta_1                  &
\alpha_1\beta_2+\alpha_2\beta_1                       &
\alpha_1\beta_3+\alpha_3\beta_1                       &
\alpha_2\beta_3+\alpha_3\beta_2                       &
\alpha_2\beta_2 \cr
b_2b_1                                                    &
\alpha_1\beta_1+\alpha_1\beta_2                  &
\alpha_1\beta_2+\alpha_2\beta_1                       &
\alpha_1\beta_3+\alpha_3\beta_1                       &
\alpha_2\beta_3+\alpha_3\beta_2                       &
\alpha_2\beta_2 \cr
b_1b_3                                                    &
\alpha_1\gamma_1+\alpha_2\gamma_1 &
\alpha_1\gamma_2+\alpha_2\gamma_1 &
\alpha_1\gamma_3+\alpha_3\gamma_1                     &
\alpha_2\gamma_3+\alpha_3\gamma_2                     &
\alpha_2\gamma_2 \cr
b_3b_1                                                    &
\alpha_1\gamma_1+\alpha_1\gamma_2 &
\alpha_1\gamma_2+\alpha_2\gamma_1 &
\alpha_1\gamma_3+\alpha_3\gamma_1                     &
\alpha_2\gamma_3+\alpha_3\gamma_2                     &
\alpha_2\gamma_2 \cr
b_2b_3                                                    &
\beta_1\gamma_1+\beta_2\gamma_1                  &
\beta_1\gamma_2+\beta_2\gamma_1                       &
\beta_1\gamma_3+\beta_3\gamma_1                       &
\beta_2\gamma_3+\beta_3\gamma_2                       &
\beta_2\gamma_2 \cr
b_3b_2                                                    &
\beta_1\gamma_1+\beta_1\gamma_2           &
\beta_1\gamma_2+\beta_2\gamma_1 & \beta_1\gamma_3+\beta_3\gamma_1
& \beta_2\gamma_3+\beta_3\gamma_2                       &
\beta_2\gamma_2 \cr
b_1^2                                                   &
\alpha_1^2+\alpha_2\alpha_3 & 0 & 0 & 0 & \alpha_2^2 \cr
b_2^2                                                   &
\beta_1^2+\beta_2\beta_3 & 0 & 0 & 0 &\beta_2^2 \cr
b_3^2                                                   &
\gamma_1^2+\gamma_2\gamma_3 & 0 & 0 & 0 & \gamma_2^2 \cr \noalign
{\hrule} }}}
$$

Since $\Delta\ne 0$, it is easy to see that the first six lines
not equal neither zero nor the last three lines. Note that the
dimension of $A^2(KG)/A^3(KG)$ equal to $5$ and $KG$ is not a
commutative algebra. From the fact $b_i^2 \equiv b_j^2\equiv 0
\pmod{A^3(KG)}$, $i\ne j$ it implies that $b_i$ linearly depends
on $b_j$, so we shall consider two interesting cases.

In the first case $b_1^2\equiv 0$, $b_2^2\equiv b_3^2\not\equiv 0
\pmod{A^3(KG)}$ and we get that $b_1\equiv \alpha_3(1+c)
\pmod{A^2(KG)}$ and by property (II) of the filtered
multiplicative $K$-basis, $b_2^2=b_3^2$. From the condition
$b_2^2\equiv b_3^2 \pmod{A^3(KG)}$ we have that
$\beta_2=\gamma_2\ne 0$ and
$(\beta_1+\gamma_1)(\beta_1+\gamma_1+\gamma_2)=0$. Since $\Delta
\ne 0$ so $\beta_1=\gamma_1+\gamma_2$ and we conclude that
$b_2=(\lambda+1)(1+a)+(1+b)+\mu(1+c)$ and
$b_3=\lambda(1+a)+(1+b)+\eta(1+c)$, where
$\lambda=\frac{\gamma_1}{\gamma_2}$,
$\mu=\frac{\beta_3}{\gamma_2}$ and
$\eta=\frac{\gamma_3}{\gamma_2}$. The fact $b_2^2=b_3^2$ gives
that $1+a^2+ab+a^3b=0$, which is impossible.

In the second case  $b_1^2\equiv b_2^2\equiv b_3^2\not\equiv 0
\pmod{A^3(KG)}$ and we can assume that $b_1b_2\equiv b_2b_1$,
$b_1b_3\equiv b_3b_1 \pmod{A^3(KG)}$ and $b_3b_2\not\equiv b_2b_3
\pmod{A^3(KG)} $. Since $b_1b_2\equiv b_2b_1$ and $b_1b_3\equiv
b_3b_1 \pmod{A^3(KG)}$ we have that
$\alpha_2\beta_1=\alpha_1\beta_2$ and
$\alpha_2\gamma_1=\alpha_1\gamma_2$. From the fact that
$b_1^2\equiv b_2^2\equiv b_3^2\not\equiv 0 \pmod{A^3(KG)}$ the
sixth column asserts that $\alpha_2=\beta_2=\gamma_2$ and the
second  column give that $\alpha_1=\beta_1=\gamma_1$, so we
conclude that $b_3b_2\equiv b_2b_3 \pmod{A^3(KG)}$, which is a
contradiction. These facts give that $KG$ has no filtered
multiplicative basis.


{\bf Case 2.} Let $G$ be the group
$$G_{24}=\gp{\quad a,b,c \quad  \vert \quad a^4=b^4=c^2=1,\;
(a,b)=(a,c)=1,\; (b,c)=a^2\quad }.$$

Using the identity
$$
      (1+c)(1+b) \equiv (1+b)(1+c)+(1+a)^2 \quad \pmod{A^3(KG)},
$$
let us compute $b_{i_1}b_{i_2}$  modulo ${A^3(KG)}$ where
($i_k=1,2,3$). The results of our  computation will be written  in
a  table, consisting of the coefficients of the decomposition
$b_{i_1}b_{i_2}$ with respect to the  basis
$$
\split \{\quad  (1+a)^{j_1}(1+b)^{j_2}(1+c)^{j_3} \quad
\mid \quad  &j_1+j_2+j_3=2;\\
& j_1,j_2=0,1,2;\; j_3=0,1 \quad \}
\endsplit
$$
of the ideal $A^2(KG)$.

$$
\eightpoint \vbox {{ \offinterlineskip
 \halign {
$#$ \; &  \vrule   \; $#$ \; &  \vrule   \; $#$ \; &  \vrule   \;
$#$ \; &  \vrule   \; $#$ \; &  \vrule   \; $#$ \; \cr
    &  (1+a)^2    &(1+a)(1+b) &  (1+a)(1+c) & (1+b)(1+c) & (1+b)^2 \cr
\noalign {\hrule}
b_1b_2                                                    &
\alpha_1\beta_1+\alpha_3\beta_2                  &
\alpha_1\beta_2+\alpha_2\beta_1                       &
\alpha_1\beta_3+\alpha_3\beta_1                       &
\alpha_2\beta_3+\alpha_3\beta_2                       &
\alpha_2\beta_2 \cr
b_2b_1                                                    &
\alpha_1\beta_1+\alpha_2\beta_3                  &
\alpha_1\beta_2+\alpha_2\beta_1                       &
\alpha_1\beta_3+\alpha_3\beta_1                       &
\alpha_2\beta_3+\alpha_3\beta_2                       &
\alpha_2\beta_2 \cr
b_1b_3                                                    &
\alpha_1\gamma_1+\alpha_3\gamma_2 &
\alpha_1\gamma_2+\alpha_2\gamma_1 &
\alpha_1\gamma_3+\alpha_3\gamma_1                     &
\alpha_2\gamma_3+\alpha_3\gamma_2                     &
\alpha_2\gamma_2 \cr
b_3b_1                                                    &
\alpha_1\gamma_1+\alpha_2\gamma_3 &
\alpha_1\gamma_2+\alpha_2\gamma_1 &
\alpha_1\gamma_3+\alpha_3\gamma_1                     &
\alpha_2\gamma_3+\alpha_3\gamma_2                     &
\alpha_2\gamma_2 \cr
b_2b_3                                                    &
\beta_1\gamma_1+\beta_3\gamma_2                  &
\beta_1\gamma_2+\beta_2\gamma_1                       &
\beta_1\gamma_3+\beta_3\gamma_1                       &
\beta_2\gamma_3+\beta_3\gamma_2                       &
\beta_2\gamma_2 \cr
b_3b_2                                                    &
\beta_1\gamma_1+\beta_2\gamma_3           &
\beta_1\gamma_2+\beta_2\gamma_1 & \beta_1\gamma_3+\beta_3\gamma_1
& \beta_2\gamma_3+\beta_3\gamma_2                       &
\beta_2\gamma_2 \cr
b_1^2                                                   &
\alpha_1^2+\alpha_2\alpha_3 & 0 & 0 & 0 & \alpha_2^2 \cr
b_2^2                                                   &
\beta_1^2+\beta_2\beta_3 & 0 & 0 & 0 &\beta_2^2 \cr
b_3^2                                                   &
\gamma_1^2+\gamma_2\gamma_3 & 0 & 0 & 0 & \gamma_2^2 \cr \noalign
{\hrule} }}}
$$
It is obvious that the first six lines not equal neither zero nor
the last three lines. Since the dimension of $A^2(KG)/A^3(KG)$
equals 5 and $KG$ is not commutative we have either $b_1b_2 \equiv
b_2b_1$, $b_1b_3 \not\equiv b_3b_1$, $b_2b_3 \not\equiv b_3b_2
\pmod{A^3(KG)}$ or $b_1b_2 \equiv b_2b_1$, $b_1b_3 \equiv b_3b_1$,
$b_2b_3 \not\equiv b_3b_2 \pmod{A^3(KG)}$, because the other cases
are analogous to these.

In the first case we get that $b_1^2\equiv b_2^2\equiv b_3^2\equiv
0 \pmod{A^3(KG)}$, so $\Delta=0$ which is impossible. In the
second case consider the following subcases:
\item{a)} $b_1^2\equiv b_2^2\equiv b_3^2\not\equiv 0 \pmod{A^3(KG)}$;
\item{b)} $b_i^2\equiv b_j^2\equiv 0$ and $b_k^2\not\equiv 0 \pmod{A^3(KG)}$;
\item{c)} $b_i^2\equiv b_j^2\not\equiv 0$ and $b_k^2\equiv 0 \pmod{A^3(KG)}$.

Since $\Delta\ne 0$ the subcase $a)$ is impossible. Consider the
subcase $b)$, and for example put $b_1^2\equiv b_2^2\equiv 0$ and
$b_3^2\not\equiv 0 \pmod{A^3(KG)}$. We get that $\alpha_2=\beta_2=
0$ and $\alpha_1=\beta_1= 0$ by second and sixth columns, so
$\Delta= 0$ which is a contradiction. The other cases also lead to
contradictions.

Assume that $b_i^2\equiv b_j^2\not\equiv 0$ and $b_k^2\equiv 0
\pmod{A^3(KG)}$, for instance $b_1^2\equiv b_2^2\not\equiv 0$ and
$b_3^2\equiv 0 \pmod{A^3(KG)}$. According to second and sixth
columns $\alpha_2=\beta_2\ne 0$ and
$(\alpha_1+\beta_1)^2=\alpha_2(\alpha_3+\beta_3)$. Since
$b_1b_2\equiv b_2b_1$ the second column gives that
$\alpha_3\beta_2=\alpha_2\beta_3$, so $\Delta= 0$ which is a
contradiction. Thus $KG$ has no a filtered multiplicative basis.

{\bf Case 3.} Let
$$
\split G=G_{27}=&\gp{\quad a,b,c \quad  \vert \quad
a^2=b^2=c^2=1,\; (a,c)=d,\; (b,c)=e,\;\\
&(a,b)=(a,d)=(a,e)=(b,d)=(b,e)=(c,d)=(c,e)=(d,e)=1 \quad }.\\
\endsplit
$$

Since
$$
M_1(G)=G, \quad M_2(G)=\gp{d,e},\quad M_3(G)=\gp{1}
$$
we obtained that $\mu(d)=\mu(e)=2$. Let us compute
$b_{i_1}b_{i_2}$ modulo ${A^3(KG)}$ where ($i_k=1,2,3$). The
results of our computation will be also written  in a table,
consisting of the coefficients of the decomposition
$b_{i_1}b_{i_2}$ with respect to the  basis
$$
\split \{\quad
(1+a)^{j_1}(1+b)^{j_2}(1+c)^{j_3}(1+d)^{j_4}(1+e)^{j_5} \quad
\mid &\quad  j_1+j_2+j_3+2j_4+2j_5=2;\\
& j_1,j_2,j_3=0,1;\; j_4,j_5=0,1\quad \}
\endsplit
$$
of the ideal $A^2(KG)$. Using the  identities:
$$
\split
(1+c)(1+a) \equiv
(1+a)(1+c)+(1+d) \pmod{A^3(KG)};\\
(1+c)(1+b) \equiv (1+b)(1+c)+(1+e) \pmod{A^3(KG)},
\endsplit
$$
we get
$$
\eightpoint
\vbox {{ \offinterlineskip
 \halign {
$#$ \;
&  \vrule   \; $#$ \;
&  \vrule   \; $#$ \;
&  \vrule   \; $#$ \;
&  \vrule   \; $#$ \;
&  \vrule   \; $#$ \;
\cr
    &  (1+a)(1+b)    &(1+a)(1+c) &  (1+b)(1+c) & (1+d) & (1+e) \cr
\noalign {\hrule}
b_1b_2                                                    &
\alpha_1\beta_2+\alpha_2\beta_1                       &
\alpha_1\beta_3+\alpha_3\beta_1                       &
\alpha_2\beta_3+\alpha_3\beta_2                       &
\alpha_3\beta_1                                       &
\alpha_3\beta_2 \cr
b_2b_1                                                    &
\alpha_1\beta_2+\alpha_2\beta_1                       &
\alpha_1\beta_3+\alpha_3\beta_1                       &
\alpha_2\beta_3+\alpha_3\beta_2                       &
\alpha_1\beta_3                                       &
\alpha_2\beta_3 \cr
b_1b_3                                                    &
\alpha_1\gamma_2+\alpha_2\gamma_1                     &
\alpha_1\gamma_3+\alpha_3\gamma_1                     &
\alpha_2\gamma_3+\alpha_3\gamma_2                     &
\alpha_3\gamma_1                                      &
\alpha_3\gamma_2 \cr
b_3b_1                                                    &
\alpha_1\gamma_2+\alpha_2\gamma_1                     &
\alpha_1\gamma_3+\alpha_3\gamma_1                     &
\alpha_2\gamma_3+\alpha_3\gamma_2                     &
\alpha_1\gamma_3                                      &
\alpha_2\gamma_3 \cr
b_2b_3                                                    &
\beta_1\gamma_2+\beta_2\gamma_1                       &
\beta_1\gamma_3+\beta_3\gamma_1                       &
\beta_2\gamma_3+\beta_3\gamma_2                       &
\beta_3\gamma_1                                       &
\beta_3\gamma_2 \cr
b_3b_2                                                    &
\beta_1\gamma_2+\beta_2\gamma_1                       &
\beta_1\gamma_3+\beta_3\gamma_1                       &
\beta_2\gamma_3+\beta_3\gamma_2                       &
\beta_1\gamma_3                                       &
\beta_2\gamma_3 \cr
b_1^2                                                   & 0 & 0 &
0 & \alpha_1\alpha_3                                      &
\alpha_2\alpha_3 \cr
b_2^2                                                   & 0 & 0 &
0 & \beta_1\beta_3                                        &
\beta_2\beta_3 \cr
b_3^2                                                   & 0 & 0 &
0 & \gamma_1\gamma_3                                      &
\gamma_2\gamma_3 \cr \noalign {\hrule} }}}
$$

It is easy to see that the first six lines not equal neither zero
nor the last three lines. Since the dimension of $A^2(KG)/A^3(KG)$
equals 5 and $KG$ is not commutative we have either $b_1b_2 \equiv
b_2b_1$, $b_1b_3 \not\equiv b_3b_1$, $b_2b_3 \not\equiv b_3b_2
\pmod{A^3(KG)}$ or $b_1b_2 \equiv b_2b_1$, $b_1b_3 \equiv b_3b_1$,
$b_2b_3 \not\equiv b_3b_2 \pmod{A^3(KG)} $, because the other
cases are similar to these.

In the first case we get that $b_1^2\equiv b_2^2\equiv b_3^2\equiv
0 \pmod{A^3(KG)}$ and $\alpha_3=\beta_3=\gamma_1=\gamma_2=0$. Let
us compute $b_{i_1}b_{i_2}b_{i_3}$ modulo ${A^4(KG)}$ where
($i_k=1,2,3$). Since the dimension of $A^3(KG)/A^4(KG)$ equal to
$7$ but we have got $8$ different elements, this case is
impossible. In the second case $b_1b_2 \equiv b_2b_1$ and $b_1b_3
\equiv b_3b_1 \pmod{A^3(KG)}$. Assume that $\alpha_3=0$. Fifth and
sixth columns give that $\beta_3=\gamma_3=0$ which is impossible,
so $\alpha_3,\beta_3,\gamma_3\ne 0$. These columns gives that
$b_2\equiv\beta_3\alpha_3^{-1}b_1 \pmod{A^2(KG)}$ which is a
contradiction, therefore $KG$ has no a filtered multiplicative
basis.

{\bf Case 4.} Let $G$ be one of the following groups:
$$
\split
G_{28}=&\gp{\quad a,b,c \quad  \vert \quad  a^4=b^2=c^2=1,\; (a,c)=a^2,\; (b,c)=d, \\
                                    &\qquad \qquad \qquad \qquad \qquad \qquad \qquad \quad \;\,(a,b)=(a,d)=(b,d)=(c,d)=1\quad };\\
G_{29}=&\gp{\quad a,b,c \quad  \vert \quad  a^4=b^2=1,\; a^2=c^2,\; (a,c)=a^2,\; (b,c)=d, \\
                                   &\qquad \qquad \qquad \qquad \qquad \qquad \qquad \quad \;\, (a,b)=(a,d)=(b,d)=(c,d)=1\quad };\\
G_{30}=&\gp{\quad a,b,c \quad  \vert \quad  a^4=b^2=c^2=1,\; (a,c)=d,\; (b,c)=a^2, \\
                                      &\qquad \qquad \qquad \qquad \qquad \qquad \qquad \quad \;\,(a,b)=(a,d)=(b,d)=(c,d)=1\quad }. \\
\endsplit
$$

\noindent If $G$ is either $G_{28}$ or $G_{29}$ then we have
$$
\split
(1+c)(1+a) &\equiv
(1+a)(1+c)+(1+a)^2 \pmod{A^3(KG)};\\
(1+c)(1+b) &\equiv
(1+b)(1+c)+(1+d) \pmod{A^3(KG)}.\\
\endsplit
$$
If $G=G_{30}$ then we have
$$
\split
(1+c)(1+a) &\equiv (1+a)(1+c)+(1+d) \pmod{A^3(KG)};\\
(1+c)(1+b) &\equiv (1+b)(1+c)+(1+a)^2 \pmod{A^3(KG)}. \\
\endsplit
$$
Using the last four identities let us compute $b_{i_1}b_{i_2}$
modulo ${A^3(KG)}$ where ($i_k=1,2,3$). The results  of our
computation will be written  in a  table as above, consisting of
the coefficients of the decomposition $b_{i_1}b_{i_2}$ with
respect to the  basis
$$
\split \{\quad  (1+a)^{j_1}(1+b)^{j_2}(1+c)^{j_3}(1+d)^{j_4} \quad
\mid \quad  &j_1+j_2+j_3+2j_4=2;\\
& j_1,j_2,j_3=0,1,2;\; j_4=0,1\quad \}
\endsplit
$$
of the ideal $A^3(KG)$:
$$
\eightpoint
\vbox {{ \offinterlineskip
 \halign {
$#$ \;
&  \vrule   \; $#$ \;
&  \vrule   \; $#$ \;
&  \vrule   \; $#$ \;
&  \vrule   \; $#$ \;
&  \vrule   \; $#$ \;
\cr
    &  (1+a)^2    &(1+a)(1+b) &  (1+a)(1+c) & (1+b)(1+c) & (1+d) \cr
\noalign {\hrule}
b_1b_2                                                    &
\alpha_1\beta_1+\Delta(\alpha,\beta)                  &
\alpha_1\beta_2+\alpha_2\beta_1                       &
\alpha_1\beta_3+\alpha_3\beta_1                       &
\alpha_2\beta_3+\alpha_3\beta_2                       &
\Omega(\alpha,\beta) \cr
b_2b_1                                                    &
\alpha_1\beta_1+\Delta(\beta,\alpha)                  &
\alpha_1\beta_2+\alpha_2\beta_1                       &
\alpha_1\beta_3+\alpha_3\beta_1                       &
\alpha_2\beta_3+\alpha_3\beta_2                       &
\Omega(\beta,\alpha) \cr
b_1b_3                                                    &
\alpha_1\gamma_1+\Delta(\alpha,\gamma)                &
\alpha_1\gamma_2+\alpha_2\gamma_1                     &
\alpha_1\gamma_3+\alpha_3\gamma_1                     &
\alpha_2\gamma_3+\alpha_3\gamma_2                     &
\Omega(\alpha,\gamma) \cr
b_3b_1                                                    &
\alpha_1\gamma_1+\Delta(\gamma,\alpha)                &
\alpha_1\gamma_2+\alpha_2\gamma_1                     &
\alpha_1\gamma_3+\alpha_3\gamma_1                     &
\alpha_2\gamma_3+\alpha_3\gamma_2                     &
\Omega(\gamma,\alpha) \cr
b_2b_3                                                    &
\beta_1\gamma_1+\Delta(\beta,\gamma)                  &
\beta_1\gamma_2+\beta_2\gamma_1                       &
\beta_1\gamma_3+\beta_3\gamma_1                       &
\beta_2\gamma_3+\beta_3\gamma_2                       &
\Omega(\beta,\gamma) \cr
b_3b_2                                                    &
\beta_1\gamma_1+\Delta(\gamma,\beta)                  &
\beta_1\gamma_2+\beta_2\gamma_1                       &
\beta_1\gamma_3+\beta_3\gamma_1                       &
\beta_2\gamma_3+\beta_3\gamma_2                       &
\Omega(\gamma,\beta) \cr
b_1^2                                                   &
\alpha_1^2+\Delta(\alpha,\alpha)                      & 0 & 0 & 0
& \Omega(\alpha,\alpha) \cr
b_2^2                                                   &
\beta_1^2+\Delta(\beta,\beta)                         & 0 & 0 & 0
& \Omega(\beta,\beta) \cr
b_3^2                                                   &
\gamma_1^2+\Delta(\gamma,\gamma)                      & 0 & 0 & 0
& \Omega(\gamma,\gamma) \cr \noalign {\hrule} }}}
$$

\noindent where if $G=G_{28}$ then
$\Delta(\delta,\epsilon)=\delta_3\epsilon_1$,
$\Omega(\delta,\epsilon)=\delta_3\epsilon_2$, if $G=G_{29}$ then
$\Delta(\delta,\epsilon)=\delta_3\epsilon_1+\delta_3\epsilon_3$,
$\Omega(\delta,\epsilon)=\delta_3\epsilon_2$
 and if $G=G_{30}$ then $\Delta(\delta,\epsilon)=\delta_3\epsilon_2$,
$\Omega(\delta,\epsilon)=\delta_3\epsilon_1$.

It is clearly that the first six lines not equal neither zero nor
the last three lines. Since the dimension of $A^2(KG)/A^3(KG)$
equal to 5 and $KG$ is not commutative we have either $b_1b_2
\equiv b_2b_1$, $b_1b_3 \not\equiv b_3b_1$, $b_2b_3 \not\equiv
b_3b_2 \pmod{A^3(KG)}$ or $b_1b_2 \equiv b_2b_1$, $b_1b_3 \equiv
b_3b_1$, $b_2b_3 \not\equiv b_3b_2 \pmod{A^3(KG)}$, because the
other cases are analogous to these.

In the first case we get that $b_1^2\equiv b_2^2\equiv b_3^2\equiv
0 \pmod{A^3(KG)}$, so $\Delta=0$ which is impossible. In the
second case $b_1b_2 \equiv b_2b_1$ and $b_1b_3 \equiv b_3b_1
\pmod{A^3(KG)}$. Assume that $\alpha_3=0$. Second and sixth
columns give that $\beta_3=\gamma_3=0$ which is impossible, so
$\alpha_3,\beta_3,\gamma_3\ne 0$. Consequences of columns $2$ and
$6$ are that $b_2\equiv\beta_3\alpha_3^{-1}b_1 \pmod{A^2(KG)}$
which is a contradiction. Thus these group algebras have no
filtered multiplicative bases.


{\bf Case 5.} Let $G$ one of the following groups:
$$
\split
G_{31}=\gp{\quad a,b,c \quad  &\vert \quad  a^4=b^4=c^2=1,\; (b,c)=a^2b^2,\; (a,c)=a^2,\; (a,b)=1\quad }; \\
G_{32}=\gp{\quad a,b,c \quad  &\vert \quad  a^4=b^4=1,\; c^2=a^2b^2,\; (b,c)=a^2b^2,\; (a,c)=a^2,\; (a,b)=1\; }; \\
G_{33}=\gp{\quad a,b,c \quad  &\vert \quad  a^4=b^4=c^2=1,\; (b,c)=a^2,\; (a,c)=a^2b^2,\; (a,b)=1\quad }; \\
G_{34}=\gp{\quad a,b,c \quad  &\vert \quad  a^4=b^4=c^2=1,\; (b,c)=b^2,\; (a,c)=a^2,\; (a,b)=1\quad };\\
G_{35}=\gp{\quad a,b,c \quad  &\vert \quad  a^4=b^4=1,\; c^2=a^2,\; (b,c)=b^2,\; (a,c)=a^2,\; (a,b)=1\quad }. \\
\endsplit
$$

Let us compute $b_{i_1}b_{i_2}$  modulo ${A^3(KG)}$,
($i_k=1,2,3$). The results  of our  computations will be written
in a  table, as before, with respect to the basis
$$
\split \{\quad  (1+a)^{j_1}(1+b)^{j_2}(1+c)^{j_3} \quad
\mid \quad  &j_1+j_2+j_3=2;\\
& j_1,j_2,j_3=0,1,2\quad \}
\endsplit
$$
of the ideal $A^3(KG)$. If $G=G_{31}$ then
$$
\split
(1+c)(1+a) &\equiv (1+a)(1+c)+(1+a)^2 \pmod{A^3(KG)}; \\
(1+c)(1+b) &\equiv (1+b)(1+c)+(1+a)^2+(1+b)^2 \pmod{A^3(KG)}, \\
\endsplit
$$
if $G=G_{32}$ then
$$
\split
(1+c)(1+a) &\equiv (1+a)(1+c)+(1+a)^2 \pmod{A^3(KG)}; \\
(1+c)(1+b) &\equiv (1+b)(1+c)+(1+a)^2+(1+b)^2 \pmod{A^3(KG)}; \\
(1+c)^2 &\equiv (1+a)^2+(1+b)^2 \pmod{A^3(KG)}, \\
\endsplit
$$
if $G=G_{33}$ then
$$
\split
(1+c)(1+a) &\equiv (1+a)(1+c)+(1+a)^2+(1+b)^2 \pmod{A^3(KG)}; \\
(1+c)(1+b) &\equiv (1+b)(1+c)+(1+a)^2 \pmod{A^3(KG)}, \\
\endsplit
$$
if $G=G_{34}$ then
$$
\split
(1+c)(1+a) &\equiv (1+a)(1+c)+(1+a)^2 \pmod{A^3(KG)}; \\
(1+c)(1+b) &\equiv (1+b)(1+c)+(1+b)^2 \pmod{A^3(KG)}, \\
\endsplit
$$
if $G=G_{35}$ then
$$
\split
(1+c)(1+a) &\equiv (1+a)(1+c)+(1+a)^2 \pmod{A^3(KG)}; \\
(1+c)(1+b) &\equiv (1+b)(1+c)+(1+b)^2 \pmod{A^3(KG)}; \\
(1+c)^2 &\equiv (1+a)^2.
\endsplit
$$

Using the last $12$ identities we get
$$
\eightpoint
\vbox {{ \offinterlineskip
 \halign {
$#$ \;
&  \vrule   \; $#$ \;
&  \vrule   \; $#$ \;
&  \vrule   \; $#$ \;
&  \vrule   \; $#$ \;
&  \vrule   \; $#$ \;
\cr
    &  (1+a)^2    &(1+a)(1+b) &  (1+a)(1+c) & (1+b)(1+c) & (1+b)^2 \cr
\noalign {\hrule}
b_1b_2                                                    &
\alpha_1\beta_1+\Delta(\alpha,\beta)                  &
\alpha_1\beta_2+\alpha_2\beta_1                       &
\alpha_1\beta_3+\alpha_3\beta_1                       &
\alpha_2\beta_3+\alpha_3\beta_2                       &
\alpha_2\beta_2+\Omega(\alpha,\beta) \cr
b_2b_1                                                    &
\alpha_1\beta_1+\Delta(\beta,\alpha)                  &
\alpha_1\beta_2+\alpha_2\beta_1                       &
\alpha_1\beta_3+\alpha_3\beta_1                       &
\alpha_2\beta_3+\alpha_3\beta_2                       &
\alpha_2\beta_2+\Omega(\beta,\alpha) \cr
b_1b_3                                                    &
\alpha_1\gamma_1+\Delta(\gamma,\alpha)                &
\alpha_1\gamma_2+\alpha_2\gamma_1                     &
\alpha_1\gamma_3+\alpha_3\gamma_1                     &
\alpha_2\gamma_3+\alpha_3\gamma_2                     &
\alpha_2\gamma_2+\Omega(\gamma,\alpha) \cr
b_3b_1                                                    &
\alpha_1\gamma_1+\Delta(\alpha,\gamma)                &
\alpha_1\gamma_2+\alpha_2\gamma_1                     &
\alpha_1\gamma_3+\alpha_3\gamma_1                     &
\alpha_2\gamma_3+\alpha_3\gamma_2                     &
\alpha_2\gamma_2+\Omega(\alpha,\gamma) \cr
b_2b_3                                                    &
\beta_1\gamma_1+\Delta(\gamma,\beta)                  &
\beta_1\gamma_2+\beta_2\gamma_1                       &
\beta_1\gamma_3+\beta_3\gamma_1                       &
\beta_2\gamma_3+\beta_3\gamma_2                       &
\beta_2\gamma_2+\Omega(\gamma,\beta) \cr
b_3b_2                                                    &
\beta_1\gamma_1+\Delta(\beta,\gamma)                  &
\beta_1\gamma_2+\beta_2\gamma_1                       &
\beta_1\gamma_3+\beta_3\gamma_1                       &
\beta_2\gamma_3+\beta_3\gamma_2                       &
\beta_2\gamma_2+\Omega(\beta,\gamma) \cr
b_1^2                                                   &
\alpha_1^2+\Delta(\alpha,\alpha)                      & 0 & 0 & 0
& \alpha_2^2+\Omega(\alpha,\alpha) \cr
b_2^2                                                   &
\beta_1^2+\Delta(\beta,\beta)                         & 0 & 0 & 0
& \beta_2^2+\Omega(\beta,\beta) \cr
b_3^2                                                   &
\gamma_1^2+\Delta(\gamma,\gamma)                      & 0 & 0 & 0
& \gamma_2^2+\Omega(\gamma,\gamma) \cr \noalign {\hrule} }}}
$$
where $\Delta(\delta,\epsilon)$ and $\Omega(\delta,\epsilon)$ is
the following:
$$
\eightpoint \vbox {{ \offinterlineskip
 \halign {\strut \; $#$ \; & \; \vrule \; $#$ \; & \;  \vrule \; $#$ \; \cr
        &   \Delta(\delta,\epsilon)   & \Omega(\delta,\epsilon) \cr \noalign {\hrule}
    G_{31} & \delta_3\epsilon_1+\delta_3\epsilon_2 &
    \delta_3\epsilon_2\cr \noalign {\hrule}
    G_{32} &
    \delta_3\epsilon_1+\delta_3\epsilon_2+\delta_3\epsilon_3 &
    \delta_3\epsilon_2+\delta_3\epsilon_3 \cr \noalign {\hrule}
     G_{33} & \delta_3\epsilon_1+\delta_3\epsilon_2 &
     \delta_3\epsilon_1\cr \noalign {\hrule}
    G_{34} & \delta_3\epsilon_1 & \delta_3\epsilon_2 \cr \noalign {\hrule}
    G_{35}  & \delta_3\epsilon_1+\delta_3\epsilon_3 &
    \delta_3\epsilon_2 \cr
 \noalign {\hrule} }}}
$$

Evidently the first six lines not equal neither zero nor the last
three lines. Since the dimension of $A^2(KG)/A^3(KG)$ equal to 5
and $KG$ is not commutative, we have either $b_1b_2 \equiv
b_2b_1$, $b_1b_3 \not\equiv b_3b_1$, $b_2b_3 \not\equiv b_3b_2
\pmod{A^3(KG)}$ or $b_1b_2 \equiv b_2b_1$, $b_1b_3 \equiv b_3b_1$,
$b_2b_3 \not\equiv b_3b_2 \pmod{A^3(KG)}$, because the other cases
are similar to these.

In both of cases we can see that $b_1b_2 \equiv b_2b_1
\pmod{A^3(KG)}$. Assume that $\alpha_3=0$. Second and sixth
columns give that $\beta_3=\gamma_3=0$ which is impossible, so
$\alpha_3,\beta_3,\gamma_3$ are not zero. Columns $2$ and $6$
imply $b_2$ depends on $b_1$ modulo $A^2(KG)$ which is a
contradiction so these group algebras have no filtered
multiplicative bases.


{\bf Case 6.} Let $G=G_{49}$ be and put $u\equiv (1+a)+(1+c)$,
$v\equiv (1+b)+(1+d)$, $w\equiv (1+b)+(1+c)+(1+d)$ and $z\equiv
(1+a)+(1+b)+(1+c) \pmod{A^2(KG)}$.

Using the identities:
$$\split
(1+b)(1+a) &\equiv  (1+a)(1+b)+(1+a)^2 \pmod{A^3(KG)};\\
(1+d)(1+c) &\equiv (1+c)(1+d)+(1+a)^2 \pmod{A^3(KG)}; \\
(1+a)^2 &\equiv (1+b)^2 \equiv (1+c)^2 \equiv
(1+d)^2 \pmod{A^3(KG)};\\
\endsplit
$$
we get that
\itemitem{$\bullet$} $\{uv,uw,uz,zu,vw,vz,wz \}$ is a basis of $A^2(KG)/A^3(KG)$;
\itemitem{$\bullet$} $\{ uzu,uvw,vzu,wzu,vuz,uzw,vwz,zuz \}$ is a basis
for $A^3(KG)/A^4(KG)$;
\itemitem{$\bullet$} $\{ vuzu,wuzu,zuzu,vzuz,wzuz,uvwz,vwzu \}$ is a basis
of $A^4(KG)/A^5(KG)$;
\itemitem{$\bullet$} $\{ vzuz,wzuz,vwuzu,vwzuz  \}$ is a basis of
$A^5(KG)/A^6(KG)$,

and the element $vwzuzu$ is a basis for $A^6(KG)$.


{\bf Case 7.} Let
$$\split G=G_{50}=\gp{\quad a,b,c,d \quad \vert \quad
a^4=b^2=&c^2=d^4=1,\,\, a^2=d^2,\\
 (a,d)=(b,c)=(c,d)=&a^2,\,\, (a,b)=(a,c)=(b,d)=1 \quad } \\
\endsplit.
$$

Using the identities:
$$
\split (1+d)(1+a) &\equiv (1+a)(1+d)+(1+a)^2 \pmod{A^3(KG)};\\
(1+c)(1+b) &\equiv (1+b)(1+c)+(1+a)^2 \pmod{A^3(KG)};\\
(1+d)(1+c) &\equiv (1+c)(1+d)+(1+a)^2 \pmod{A^3(KG)};\\
(1+a)^2 &\equiv (1+d)^2 \pmod{A^3(KG)},\\
\endsplit
$$
let us compute $b_{i_1}b_{i_2}$ modulo ${A^3(KG)}$ where
$i_k=1,2,3,4 \}$. The results  of our computations we shall write
in a table, similar to previous cases with respect to the basis
$$
\split \{\quad  (1+a)^{j_1}(1+b)^{j_2}(1+c)^{j_3}(1+d)^{j_4} \quad
\mid \quad  &j_1+j_2+j_3+j_4=2;\\
& j_1=0,1,2;\; j_2,j_3,j_4=0,1\quad \}
\endsplit
$$
of the ideal $A^2(KG)$:

$$
\eightpoint \vbox {{ \offinterlineskip
 \halign {
$#$ &  \vrule       $#$  &  \vrule       $#$  & \vrule  $#$  &
\vrule       $#$  & \vrule    $#$  & \vrule    $#$ & \vrule  $#$ &
\vrule    $#$  \cr
 &
\scriptstyle{(1+a)(1+b)} &\scriptstyle{(1+a)(1+c)} &
\scriptstyle{(1+a)(1+d)} & \scriptstyle{(1+b)(1+c)}
&\scriptstyle{(1+b)(1+d)}& \scriptstyle{(1+c)(1+d)} &
\scriptstyle{(1+a)^2} \cr \noalign {\hrule}
\scriptstyle{ b_1b_2 }& \scriptstyle{\Delta^{1,2}(\alpha,\beta)} &
\scriptstyle{\Delta^{1,3}(\alpha,\beta)} &
\scriptstyle{\Delta^{1,4}(\alpha,\beta)} &
\scriptstyle{\Delta^{2,3}(\alpha,\beta)} &
\scriptstyle{\Delta^{2,4}(\alpha,\beta)} &
\scriptstyle{\Delta^{3,4}(\alpha,\beta)} &
\scriptstyle{\Omega_{\alpha,\beta}+\alpha_3\beta_2+\alpha_4\beta_1+\alpha_4\beta_3}
\cr
\scriptstyle{ b_2b_1 }& \scriptstyle{\Delta^{1,2}(\alpha,\beta)} &
\scriptstyle{\Delta^{1,3}(\alpha,\beta)} &
\scriptstyle{\Delta^{1,4}(\alpha,\beta)} &
\scriptstyle{\Delta^{2,3}(\alpha,\beta)} &
\scriptstyle{\Delta^{2,4}(\alpha,\beta)} &
\scriptstyle{\Delta^{3,4}(\alpha,\beta)} &
\scriptstyle{\Omega_{\alpha,\beta}+\alpha_2\beta_3+\alpha_1\beta_4+\alpha_3\beta_4}
\cr
\scriptstyle{ b_1b_3 }& \scriptstyle{\Delta^{1,2}(\alpha,\gamma)}
& \scriptstyle{\Delta^{1,3}(\alpha,\gamma)} &
\scriptstyle{\Delta^{1,4}(\alpha,\gamma)} &
\scriptstyle{\Delta^{2,3}(\alpha,\gamma)} &
\scriptstyle{\Delta^{2,4}(\alpha,\gamma)} &
\scriptstyle{\Delta^{3,4}(\alpha,\gamma)} &
\scriptstyle{\Omega_{\alpha,\gamma}+\alpha_3\gamma_2+\alpha_4\gamma_1+\alpha_4\gamma_3}
\cr
\scriptstyle{ b_3b_1 }& \scriptstyle{\Delta^{1,2}(\alpha,\gamma)}
& \scriptstyle{\Delta^{1,3}(\alpha,\gamma)} &
\scriptstyle{\Delta^{1,4}(\alpha,\gamma)} &
\scriptstyle{\Delta^{2,3}(\alpha,\gamma)} &
\scriptstyle{\Delta^{2,4}(\alpha,\gamma)} &
\scriptstyle{\Delta^{3,4}(\alpha,\gamma)} &
\scriptstyle{\Omega_{\alpha,\gamma}+\alpha_2\gamma_3+\alpha_1\gamma_4+\alpha_3\gamma_4}
\cr
\scriptstyle{ b_1b_4 }& \scriptstyle{\Delta^{1,2}(\alpha,\delta)}
& \scriptstyle{\Delta^{1,3}(\alpha,\delta)} &
\scriptstyle{\Delta^{1,4}(\alpha,\delta)} &
\scriptstyle{\Delta^{2,3}(\alpha,\delta)} &
\scriptstyle{\Delta^{2,4}(\alpha,\delta)} &
\scriptstyle{\Delta^{3,4}(\alpha,\delta)} &
\scriptstyle{\Omega_{\alpha,\delta}+\alpha_3\delta_2+\alpha_4\delta_1+\alpha_4\delta_3}
\cr
\scriptstyle{ b_4b_1 }& \scriptstyle{\Delta^{1,2}(\alpha,\delta)}
& \scriptstyle{\Delta^{1,3}(\alpha,\delta)} &
\scriptstyle{\Delta^{1,4}(\alpha,\delta)} &
\scriptstyle{\Delta^{2,3}(\alpha,\delta)} &
\scriptstyle{\Delta^{2,4}(\alpha,\delta)} &
\scriptstyle{\Delta^{3,4}(\alpha,\delta)} &
\scriptstyle{\Omega_{\alpha,\delta}+\alpha_2\delta_3+\alpha_1\delta_4+\alpha_3\delta_4}
\cr
\scriptstyle{ b_2b_3 }& \scriptstyle{\Delta^{1,2}(\beta,\gamma)} &
\scriptstyle{\Delta^{1,3}(\beta,\gamma)} &
\scriptstyle{\Delta^{1,4}(\beta,\gamma)} &
\scriptstyle{\Delta^{2,3}(\beta,\gamma)} &
\scriptstyle{\Delta^{2,4}(\beta,\gamma)} &
\scriptstyle{\Delta^{3,4}(\beta,\gamma)} &
\scriptstyle{\Omega_{\beta,\gamma}+\beta_3\gamma_2+\beta_4\gamma_1+\beta_4\gamma_3}
\cr
\scriptstyle{ b_3b_2 }& \scriptstyle{\Delta^{1,2}(\beta,\gamma)} &
\scriptstyle{\Delta^{1,3}(\beta,\gamma)} &
\scriptstyle{\Delta^{1,4}(\beta,\gamma)} &
\scriptstyle{\Delta^{2,3}(\beta,\gamma)} &
\scriptstyle{\Delta^{2,4}(\beta,\gamma)} &
\scriptstyle{\Delta^{3,4}(\beta,\gamma)} &
\scriptstyle{\Omega_{\beta,\gamma}+\beta_2\gamma_3+\beta_1\gamma_4+\beta_3\gamma_4}
\cr
\scriptstyle{ b_2b_4 }& \scriptstyle{\Delta^{1,2}(\beta,\delta)} &
\scriptstyle{\Delta^{1,3}(\beta,\delta)} &
\scriptstyle{\Delta^{1,4}(\beta,\delta)} &
\scriptstyle{\Delta^{2,3}(\beta,\delta)} &
\scriptstyle{\Delta^{2,4}(\beta,\delta)} &
\scriptstyle{\Delta^{3,4}(\beta,\delta)} &
\scriptstyle{\Omega_{\beta,\delta}+\beta_3\delta_2+\beta_4\delta_1+\beta_4\delta_3}
\cr
\scriptstyle{ b_4b_2 }& \scriptstyle{\Delta^{1,2}(\beta,\delta)} &
\scriptstyle{\Delta^{1,3}(\beta,\delta)} &
\scriptstyle{\Delta^{1,4}(\beta,\delta)} &
\scriptstyle{\Delta^{2,3}(\beta,\delta)} &
\scriptstyle{\Delta^{2,4}(\beta,\delta)} &
\scriptstyle{\Delta^{3,4}(\beta,\delta)} &
\scriptstyle{\Omega_{\beta,\delta}+\beta_2\delta_3+\beta_1\delta_4+\beta_3\delta_4}
\cr
\scriptstyle{ b_3b_4 }& \scriptstyle{\Delta^{1,2}(\gamma,\delta)}
& \scriptstyle{\Delta^{1,3}(\gamma,\delta)} &
\scriptstyle{\Delta^{1,4}(\gamma,\delta)} &
\scriptstyle{\Delta^{2,3}(\gamma,\delta)} &
\scriptstyle{\Delta^{2,4}(\gamma,\delta)} &
\scriptstyle{\Delta^{3,4}(\gamma,\delta)} &
\scriptstyle{\Omega_{\gamma,\delta}+\gamma_3\delta_2+\gamma_4\delta_1+\gamma_4\delta_3}
\cr
\scriptstyle{ b_4b_3 }& \scriptstyle{\Delta^{1,2}(\gamma,\delta)}
& \scriptstyle{\Delta^{1,3}(\gamma,\delta)} &
\scriptstyle{\Delta^{1,4}(\gamma,\delta)} &
\scriptstyle{\Delta^{2,3}(\gamma,\delta)} &
\scriptstyle{\Delta^{2,4}(\gamma,\delta)} &
\scriptstyle{\Delta^{3,4}(\gamma,\delta)} &
\scriptstyle{\Omega_{\gamma,\delta}+\gamma_2\delta_3+\gamma_1\delta_4+\gamma_3\delta_4}
\cr
\scriptstyle{ b_1^2 }& \scriptstyle{0} & \scriptstyle{0} &
\scriptstyle{0} & \scriptstyle{0} & \scriptstyle{0} &
\scriptstyle{0} &
\scriptstyle{\Omega_{\alpha,\alpha}+\alpha_2\alpha_3+\alpha_1\alpha_4+\alpha_3\alpha_4}
\cr
\scriptstyle{ b_2^2 }& \scriptstyle{0} & \scriptstyle{0} &
\scriptstyle{0} & \scriptstyle{0} & \scriptstyle{0} &
\scriptstyle{0} &
\scriptstyle{\Omega_{\beta,\beta}+\beta_2\beta_3+\beta_1\beta_4+\beta_3\beta_4}
\cr
\scriptstyle{ b_3^2 }& \scriptstyle{0} & \scriptstyle{0} &
\scriptstyle{0} & \scriptstyle{0} & \scriptstyle{0} &
\scriptstyle{0} &
\scriptstyle{\Omega_{\gamma,\gamma}+\gamma_2\gamma_3+\gamma_1\gamma_4+\gamma_3\gamma_4}
\cr
\scriptstyle{ b_4^2 }& \scriptstyle{0} & \scriptstyle{0} &
\scriptstyle{0} & \scriptstyle{0} & \scriptstyle{0} &
\scriptstyle{0} &
\scriptstyle{\Omega_{\delta,\delta}+\delta_2\delta_3+\delta_1\delta_4+\delta_3\delta_4}
\cr \noalign {\hrule} }}}
$$
where $\Omega_{\alpha,\beta}=\alpha_1\beta_1+\alpha_4\beta_4$ and
$\Delta^{i,j}(\alpha,\beta)=\alpha_i\beta_j+\alpha_j\beta_i$.

It is easy to see that the first twelve lines not equal neither
zero nor the last four lines.

Since $\Delta^{i,j}(\varepsilon,\eta)$ is a subdeterminant of
$\Delta$ and $\Delta\ne 0$, by expansion theorem of determinant
$b_ib_j$ cannot be equivalent other else $ b_kb_l \pmod{A^3(KG)}$
apart from the case when $k=j$ and $l=i$.

Assume that $\{ u\equiv b_1,v\equiv b_2,w\equiv b_3,z\equiv b_4
\pmod{A^2(KG)} \}$ and the coefficients of $u,v,w,z$ will be
denoted by $\alpha_i, \beta_i, \gamma_i, \delta_i$, respectively.


Since the dimension of $A^2(KG)/A^3(KG)$ is equal to 7 and this
group algebra is not commutative, we have that
$$
\aligned uv \equiv vu, \quad uw &\equiv wu, \quad uz\equiv zu,
\quad vw\equiv wv, \quad vz\equiv
zv, \\
wz &\not\equiv zw, \quad u^2\equiv v^2\equiv w^2\equiv z^2 \equiv
0 \pmod{A^3(KG)},
\endaligned \tag 6
$$
and the other cases are analogous to this one.

Assume that $KG$ has a filtered multiplicative basis and  $\{
u,v,w,z \}$ is a basis of $A(KG)/A^2(KG)$, satisfies $(6)$ and
$$
\aligned &\{ uv,uw,uz,vw,vz,wz,zw \} \quad \text{is a basis
for} \quad A^2(KG)/A^3(KG);\\
&\{ uvw,uvz,uwz,uzw,vwz,vzw,wzw,zwz \} \quad \text{is a basis
of} \quad A^3(KG)/A^4(KG);\\
&\{ uvwz,uvzw,uwzw,uzwz,vwzw,vzwz,wzwz \} \, \text{is a basis
of} \,  A^4(KG)/A^5(KG);\\
&\{ uvwzw,uvzwz,uwzwz,vwzwz \} \quad \text{is a basis
for} \quad A^5(KG)/A^6(KG);\\
&\{ uvwzwz \} \quad \text{is a basis
for} \quad A^6(KG).\\
\endaligned \tag 7
$$

Suppose that $\alpha_4=0$ and there exists $b\in \{v,w,z \}$ such
that $b$ is congruent with
$\varepsilon_1(1+a)+\varepsilon_2(1+b)+\varepsilon_3(1+c)+\varepsilon_4(1+d)
\pmod{A^2(KG)}$ and $\varepsilon_4=0$. The facts $u^2\equiv
b^2\equiv 0$ and $ub\equiv bu \pmod{A^3(KG)}$ give that
$\alpha_3\varepsilon_2+\alpha_2\varepsilon_3=0$ and
$\alpha_1^2+\alpha_2\alpha_3=\varepsilon_1^2+\varepsilon_2\varepsilon_3=0$.
It is very simple to prove that either $b_1\equiv 0$ or $b_1\equiv
b_i \pmod{A^2(KG)}$, which is impossible.

Now, we shall consider two subcases.

{\bf Subcase $1$.} Suppose that $\alpha_4=0$ and
$\beta_4=\gamma_4=1$. For $\alpha_2=0$ it follows that $\Delta=0$,
so we can also assume that $\alpha_2=1$. According to eighth
column of the previous table
$$
\aligned
\alpha_1^2+\alpha_3&=0;\\
\beta_1+\beta_3+1&=\beta_1^2+\beta_2\beta_3;\\
\gamma_1+\gamma_3+1&=\gamma_1^2+\gamma_2\gamma_3;\\
\endaligned \tag 8
$$
Since $vw\equiv wv \pmod{A^3(KG)}$ we get
$\beta_3\gamma_2+\gamma_1+\gamma_3+1=\beta_2\beta_3+\beta_1+\beta_3+1$
and using $(8)$ it follows that
$$
(\beta_1+\gamma_1)^2=(\beta_2+\gamma_2)(\beta_3+\gamma_3). \tag 9
$$
Also eighth column of the previous table and $uv\equiv vu$,
$uw\equiv wu \pmod{A^3(KG)}$ give that
$\alpha_1^2\beta_2+\beta_3=\alpha_1^2\gamma_2+\gamma_3$, so
$$
\alpha_1^2(\beta_2+\gamma_2)=\beta_3+\gamma_3. \tag 10
$$
Thus $(9)$ and $(10)$ give the equation
$\beta_1+\gamma_1=\alpha_1(\beta_2+\gamma_2)$.

Since $\alpha_1^2=\alpha_3$ we have established $v+w\equiv
(\beta_2+\gamma_2)u \pmod{A^3(KG)}$ which is a contradiction.

{\bf Subcase $2$.} Suppose that $\alpha_4,\beta_4\ne 0$ and
without loss of generality we can assume that
$\alpha_4=\beta_4=1$. Simple computations show that $\{ u\equiv
b_1+b_2,v\equiv b_2,w\equiv b_3,z\equiv b_4 \pmod{A^2(KG)} \}$
form a basis of $A(KG)/A^2(KG)$, satisfies conditions $(6)$ and
$(7)$, but it is a contradiction to subcase $1$, so this group
algebra has no filtered multiplicative basis. This completes the
proof of the theorem.


\enddocument